\documentclass{amsart}

\usepackage{amssymb}
\usepackage{graphicx}

\newcommand{\ovl}{\overline}
\newcommand{\ie}{{i.e.\ }}

\newcommand{\mA}{{\mathcal A}}
\newcommand{\mB}{{\mathcal B}}
\newcommand{\mC}{{\mathcal C}}
\newcommand{\mQ}{{\mathcal Q}}

\newcommand{\N}{\mathbb{N}}  
\newcommand{\Z}{\mathbb{Z}}  
\newcommand{\R}{\mathbb{R}}  
\newcommand{\C}{\mathbb{C}}  
\newcommand{\iC}{\mathrm{i}} 
\newcommand{\Eu}{\mathrm{e}} 

\newcommand{\les}{\lesssim}
\newcommand{\ges}{\gtrsim}

\newcommand{\sss}{\approx}

\newcommand{\1}[1]{\frac{1}{#1}}
\newcommand{\nd}{\frac{n}2}

\newcommand{\abs}[1]{\left|#1\right|}
\newcommand{\norm}[1]{\left\|#1\right\|}
\newcommand{\tonde}[1]{\left(#1\right)}
\newcommand{\quadre}[1]{\left[#1\right]}
\newcommand{\graffe}[1]{\left\{#1\right\}}

\newcommand{\changeto}{\longrightarrow}

\newcommand{\de}{\partial}
\newcommand{\der}[3][]%
{{\frac{{\rm d}^{#1} #2}{{\rm d}{#3}^{#1}}}}
\newcommand{\pder}[3][]%
{{\frac{\de^{#1} #2}{\de{#3}^{#1}}}}
\renewcommand{\d}{\,{\rm d}} 

\newcommand{\qt}{{\widetilde{q}}}
\newcommand{\rt}{{\widetilde{r}}}

\newcommand{\union}{\cup}
\newcommand{\compose}{\circ}

\DeclareMathOperator{\esssup}{ess\,sup}

\theoremstyle{plain}
\newtheorem{theorem}{Theorem}[section]
\newtheorem{lemma}[theorem]{Lemma}
\newtheorem{proposition}[theorem]{Proposition}
\theoremstyle{definition}
\newtheorem{example}[theorem]{Example}
\theoremstyle{remark}
\newtheorem{remark}[theorem]{Remark}

\begin{document}

\title[Trigonometric sums and oscillatory integrals]
{Some remarks on the $L^p-L^q$ boundedness
  of trigonometric sums and oscillatory integrals}
\author{Damiano Foschi}
\address{%
  Dipartimento di Matematica Pura e Applicata \\
  Universit\`a di L'Aquila}
\email{foschi@univaq.it}
\urladdr{http://univaq.it/\textasciitilde foschi/}
\date{October 15, 2003}


\begin{abstract}
  We discuss the asymptotic behaviour
  for the best constant
  in $L^p$-$L^q$ estimates for trigonometric polinomials
  and for an integral operator which is related
  to the solution of inhomogeneous Schr\"odinger equations.
  This gives us an opportunity to review 
  some basic facts about oscillatory integrals
  and the method of stationary phase,
  and also to make some remarks in connection with
  Strichartz estimates.
\end{abstract}

\maketitle

\section{Introduction}

Let $u(t,x)$ be the solution
of the homogeneous Schr\"odinger equation
\begin{equation*}
  \iC \de_t u - \Delta u = 0,
\end{equation*}
with initial data $u(0,x)=f(x)$.
Let $v(t,x)$ be the solution
of the inhomogeneous Schr\"odinger equation
\begin{equation*}
  \iC \de_t v - \Delta v = F(t,x),
\end{equation*}
with zero initial data.
It is known~\cite{Obe1989, Har1990, Kat1994} that
inhomogeneous Strichartz estimates of the form
\begin{equation*}
  \norm{v}_{L^q(\R; L^r(\R^n))} \les
  \norm{F}_{L^{\qt'}(\R; L^{\rt}(\R^n))}
\end{equation*}
are valid even for some pairs of exponents $(q,r)$, $(\qt,\rt)$
which are not admissible
for the homogeneous Strichartz estimate
\begin{equation} \label{eq:1}
  \norm{u}_{L^q(\R; L^r(\R^n))} \les \norm{f}_{L^2(\R^n)}.
\end{equation}
While searching for counterexamples which could help us
understand what
the optimal range for the exponents $q,r,\qt,\rt$
in~\eqref{eq:1} could be,
a simplification of the problem led us to consider
the integral operator $T: L^p(0,1) \to L^q(0,1)$ defined by
\begin{equation} \label{eq:2}
  Tf(t) = \int_0^1 \Eu^{\iC N/(1+t+s)} f(s) \d s.
\end{equation}
A further simplified discrete version of this integral operator
is represented by the operator which assigns
to $N$ complex numbers $a_0, \dots, a_{N-1}$
the trigonometric polynomial
$\sum_{n=0}^{N-1} a_n \Eu^{\iC nt}$,
acting from $\ell^p(\C^N)$ to $L^q(-\pi,\pi)$.
We are interested in the asymptotic behaviour
of its operator norm as $N \to \infty$.
This becomes an interesting exercise
in elementary harmonic analysis
whose solution (theorem~\ref{theorem:1}) is discussed in
sections~\ref{sec:trigonometric-sums},%
~\ref{sec:upper-bounds},~\ref{sec:lower-bounds}.
In section~\ref{sec:an-integral-operator}
we obtain estimates for integral operators like~\eqref{eq:2}.
We then use them in section~\ref{sec:local-estim-inhom}
to find the optimal range of exponent for
a weaker local version of Strichartz estimates
(theorem~\ref{theorem:3}).
The details of the proofs of the various lemmata
about oscillatory integrals
which are needed throughout the paper
are collected in section~\ref{sec:proof-lemmata}.

\section{Trigonometric sums} \label{sec:trigonometric-sums}

Given $N$ complex numbers $a_0, a_2, \dots, a_{N-1}$,
the trigonometric sum
\begin{equation}
  \label{eq:3}
  f(t) = \sum_{n=0}^{N-1} a_n \Eu^{\iC n t}
\end{equation}
defines a smooth $2\pi$-periodic function.
Let $T_N$ be the linear operator
from $\C^N$ to $C^\infty\tonde{[-\pi,\pi]; \C}$
which maps the vector $a=(a_0, a_1, \dots, a_{N-1})$
to the function $f$.
For any $p,q \in [1,\infty]$, let us denote by $C_N(p \to q)$
the best constant for the estimate
\begin{equation*}
  \norm{f}_{L^q} \le C \norm{a}_{\ell^p},
\end{equation*}
more precisely
\begin{equation*}
  C_N(p \to q) = \norm{T_N}_{\ell^p \to L^q} =
  \sup_{a \in \C^N \setminus \graffe{0}}
  \frac{\norm{T_N(a)}_{L^q}}{\norm{a}_{\ell^p}},
\end{equation*}
where the norms are defined by
\begin{align*}
  \norm{a}_{\ell^p} &=
  \tonde{\sum_{n=0}^{N-1} \abs{a_n}^p}^{1/p},
  1\le p<\infty, &
  \norm{a}_{\ell^\infty} &= \max_{0\le n\le N-1} \abs{a_n},
\end{align*}
and
\begin{align*}
  \norm{f}_{L^q} &=
  \tonde{\1{2\pi} \int_{-\pi}^{\pi} \abs{f(t)}^q}^{1/q},
  1\le q<\infty, &
  \norm{f}_{L^\infty} &= \max_{\abs{t}\le\pi} \abs{f(t)}.
\end{align*}

{\bf Problem:} Is it possible to compute $C_N(p \to q)$,
or at least to describe its asymptotic behaviour
as $N \to \infty$?

\begin{remark}
  Since the operator $T_N$
  is defined on a finite dimensional vector space, 
  we know that the constant $C_N(p \to q)$ is always finite
  and that for any choice of $N, p, q$
  there esists some \emph{maximizer} $a \in \C^N$
  for which we have
  $\norm{T_N(a)}_{L^q} = C_N(p \to q) \cdot \norm{a}_{\ell^p}$.
\end{remark}

\begin{remark}
  In order to facilitate
  the visualization of relations among the various estimates,
  it will be convenient to use the notation $(p \to q)$
  to indicate the point $(1/p, 1/q)$
  in the unit square $\mQ = [0,1]^2$
  and view $C_N$ as a function defined on $\mQ$.
\end{remark}

\begin{figure}
  \centering
  \begin{picture}(0,0)%
    \includegraphics{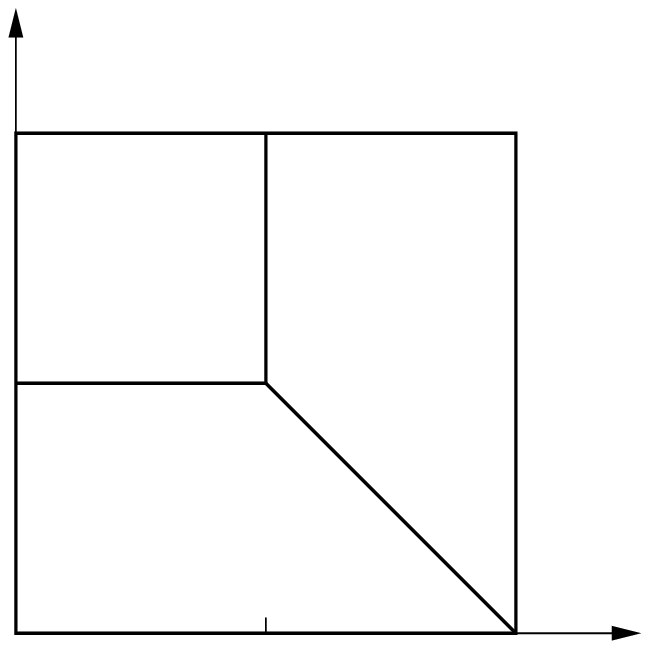}%
  \end{picture}%
  \setlength{\unitlength}{3947sp}%
  \begin{picture}(3237,3265)(376,-2414)
    \put(2079,-916){\makebox(0,0)[lb]{\smash{$C_N=1$}}}
    \put(841,-1861){\makebox(0,0)[lb]{\smash{%
          $C_N \sss N^{1-\1q-\1p}$}}}
    \put(1148,-166){\makebox(0,0)[lb]{\smash{$\mC$}}}
    \put(713,-646){\makebox(0,0)[lb]{\smash{%
          $C_N \sss N^{\12-\1p}$}}}
    \put(2281,-182){\makebox(0,0)[lb]{\smash{$\mA$}}}
    \put(1321,-1381){\makebox(0,0)[lb]{\smash{$\mB$}}}
    \put(391,-1013){\makebox(0,0)[lb]{\smash{$\12$}}}
    \put(443,209){\makebox(0,0)[lb]{\smash{$1$}}}
    \put(376,568){\makebox(0,0)[lb]{\smash{$\1q$}}}
    \put(1756,-2356){\makebox(0,0)[lb]{\smash{$\12$}}}
    \put(3233,-2356){\makebox(0,0)[lb]{\smash{$\1p$}}}
    \put(2926,-2355){\makebox(0,0)[lb]{\smash{$1$}}}
  \end{picture}
  \caption{The region
    $\mQ = \mA \union \mB \union \mC = [0,1] \times [0,1]$.}
  \label{fig:QABC}
\end{figure}

Let's decompose the square $\mQ$ into the three regions
(see figure~\ref{fig:QABC})
\begin{align*}
  \mA &= \graffe{(p \to q):
    \12 \le \1p \le 1, \quad 1-\1p \le \1q \le 1}, \\
  \mB &= \graffe{(p \to q): 
    0 \le \1q \le \12, \quad 0 \le \1p \le 1-\1q}, \\
  \mC &= \graffe{(p \to q): 
    0 \le \1p \le \12, \quad \12 \le \1q \le 1}.
\end{align*}

In sections~\ref{sec:upper-bounds} and~\ref{sec:lower-bounds}
we calculate upper and lower bounds
for $C_N$ which are summarized in the following theorem.

\begin{theorem} \label{theorem:1}
There exists positive absolute constants $c_{\mB}$ and $c_{\mC}$
such that:
  \begin{align*}
    (p \to q) \in \mA &\implies
    C_N(p \to q) = 1; \\
    (p \to q) \in \mB &\implies
    c_{\mB} N^{1-\1q-\1p} \le C_N(p \to q) \le N^{1-\1q-\1p}; \\
    (p \to q) \in \mC &\implies
    c_{\mC} N^{\12 - \1p} \le C_N(p \to q) \le N^{\12 - \1p}.
  \end{align*}
\end{theorem}

\section{Upper bounds} \label{sec:upper-bounds}

\begin{remark}[Trivial \emph{dispersive} estimate]
  If we ignore the oscillations,
  take absolute values
  and use triangular inequality in~\eqref{eq:3},
  we immediately obtain the $(1 \to \infty)$ estimate,
  $\norm{f}_{L^\infty} \le \norm{a}_{\ell^1}$.
  This implies that 
  \begin{equation}
    \label{eq:4}
    C_N(1 \to \infty) \le 1.
  \end{equation}
\end{remark}

\begin{remark}[\emph{Energy} estimate] \label{rem:energy}
  We can exploit $L^2$ orthogonality
  of the oscillating terms in~\eqref{eq:3}
  and obtain the $(2 \to 2)$ estimate,
  \begin{multline*}
    \norm{f}_{L^2}^2 =
    \1{2\pi} \int_{-\pi}^\pi f(t) \ovl{f(t)} \d t = 
    \sum_{k,h} \frac{a_k \ovl{a_h}}{2\pi}
    \int_{-\pi}^\pi \Eu^{\iC (k-h) t} \d t = \\
    = \sum_{k,h} a_k \ovl{a_h} \delta(k-h) = 
    \sum_{k} a_k \ovl{a_k} = \norm{a}_{\ell^2}^2.
  \end{multline*}
  This implies that
  \begin{equation}
    \label{eq:5}
    C_N(2 \to 2) = 1.
  \end{equation}
\end{remark}

\begin{remark}[Interpolation]
  The Riesz-Thorin interpolation theorem (\cite{BerLof1976})
  applied to our operator $T_N$ tells us that
  \begin{equation*}
    C_N(p \to q) \le
    C_N(p_0 \to q_0)^{1-\theta} C_N(p_1 \to q_1)^\theta,
  \end{equation*}
  when
  \begin{align*}
    \1p &= \frac{1-\theta}{p_0} + \frac{\theta}{p_1}, &
    \1q &= \frac{1-\theta}{q_0} + \frac{\theta}{q_1}, &
    & 0 \le \theta \le 1.
  \end{align*}
  This amounts to saying that
  $\log C_N$ is a convex function on $\mQ$.
  In particular, if we have bounds for $C_N$
  at any two points $X$ and $Y$ of $\mQ$,
  then interpolation gives us bounds for $C_N$
  on the whole segment in $\mQ$ connecting $X$ with $Y$.
\end{remark}

\begin{remark}[H\"older inclusions]
  If $1\le q \le \widetilde{q} \le \infty$ 
  we can apply H\"older's inequality to the norm of $f$,
  $\norm{f}_{L^q} \le \norm{f}_{L^{\widetilde{q}}}$,
  and obtain the following condition for $C_N$:
  \begin{equation*}
    q \le \widetilde{q} \implies
    C_N(p \to q) \le C_N(p \to \widetilde{q}).
  \end{equation*}
  If $1\le p \le \widetilde{p} \le \infty$ 
  we can apply H\"older's inequality to the norm of $a$,
  $\norm{a}_{\ell^p} \le
  N^{1/p - 1/{\widetilde{p}}} \norm{a}_{\ell^{\widetilde{p}}}$,
  and obtain the following condition for $C_N$:
  \begin{equation*}
    p \le \widetilde{p} \implies
    C_N(\widetilde{p} \to q) \le
    N^{\1p - \1{\widetilde{p}}} C_N(p \to q).
  \end{equation*}
  Looking at the $\mQ$ square,
  these conditions mean that
  an upper bound at one point in $\mQ$
  implies upper bounds at any point which can be reached
  by moving upward or leftward.
\end{remark}

Combining together the above remarks
we obtain the following upper bounds for $C_N$:

\begin{proposition}
  We have
  \begin{equation}
    \label{eq:6}
    C_N(p \to q) \le
    \max\graffe{1,\; N^{1-\1q-\1p},\; N^{\12-\1p}}.
  \end{equation}
\end{proposition}

\begin{proof}
  Interpolation between
  the $(1 \to \infty)$ estimate~\eqref{eq:4}
  and the $(2 \to 2)$ estimate~\eqref{eq:5}
  proves the result when $1/q = 1-1/p \in [0,1/2]$.
  All other cases follow from these
  by applying H\"older's inequality.
\end{proof}

\begin{remark}
  It is interesting to note that
  we had to look at the structure of the operator $T_N$
  only for the \emph{dispersive} $(1 \to \infty)$ estimate
  and the \emph{energy} $(2 \to 2)$ estimate.
  All other estimates followed from these two cases
  using only the structure and interpolation properties
  of $L^p$ spaces,
  without having to look at the structure of the operator $T_N$.
  A similar situation happens when we want to prove 
  Strichartz estimates for dispersive evolution operators
  \cite{KeeTao1998}.
\end{remark}

\section{Lower bounds} \label{sec:lower-bounds}

We can obtain lower bounds for $C_N$
by computing the norms of $f$ and $a$
for specific examples.

\subsection{Lower bounds for region $\mA$}

\begin{example}[Trivial case]
  Let us choose $a = (1,0,\dots,0)$.
  In this case $f(t) = 1$ is constant.
  Hence, for any $p$ and $q$ we have
  \begin{equation*}
    C_N(p \to q) \ge \frac{\norm{f}_{L^q}}{\norm{a}_{\ell^p}}
    = \frac{1}{1} = 1.
  \end{equation*}
  Together with the upper bound~\eqref{eq:6},
  this proves that $C_N(p \to q) = 1$ in the region $\mA$.
  Moreover, in this region any choice of $a \in \C^N$
  whose components are all vanishing except for one
  is a maximizer.
\end{example}

\subsection{Lower bounds for region $\mB$}

\begin{example}[Dirichlet kernels] \label{exa:dirichlet}
  Let us choose $a=(1,1,\dots,1)$.
  Its norm is $\norm{a}_{\ell^p} = N^{1/p}$.
  We can compute $f$ explicitly,
  \begin{equation*}
    f(t) = \sum_{n=0}^{N-1} \Eu^{\iC n t}
    = \frac{\Eu^{\iC N t}-1}{\Eu^{\iC t}-1}
    = \Eu^{\iC \frac{N-1}{2} t} D_N(t),
  \end{equation*}
  where $D_N$ is the Dirichlet kernel
  \begin{equation}
    \label{eq:7}
    D_N(t) = \frac{\sin\tonde{Nt/2}}{\sin(t/2)}.
  \end{equation}
  To estimate $D_N$ from below, 
  we use the fact that 
  \begin{equation*}
    0 \le \alpha \le \frac\pi2 \implies 
    \frac{2\alpha}\pi \le \sin \alpha \le \alpha
  \end{equation*}
  and find
  \begin{equation*}
    \abs{t} \le \frac{\pi}{N} \implies
    \abs{D_N(t)} \ge \frac{2N}{\pi}.
  \end{equation*}
  It follows that
  \begin{equation}
    \label{eq:8}
    \norm{f}_{L^q} = \norm{D_N}_{L^q} \ge 
    \frac{2N}{\pi} \tonde{\1{2\pi} \cdot \frac{2\pi}{N}}^{1/q} =
    \frac{2}{\pi} N^{1 - \1q},
  \end{equation}
  from which we obtain
  \begin{equation*}
    C_N(p \to q) \ge
    \frac{\norm{f}_{L^q}}{\norm{a}_{\ell^p}} \ge
    \frac{2}{\pi} N^{1 - \1q - \1p}.
  \end{equation*}
  In particular,
  this example shows that in region $\mB$
  the exponent $1-1/q-1/p$
  which appears in the upper bound~\eqref{eq:6}
  is sharp
  and that we can take $c_\mB = 2/\pi$
  in theorem~\ref{theorem:1}.
\end{example}

The following lemma,
which we prove in section~\ref{sec:proof-lemmata},
improves the estimate~\eqref{eq:8}.

\begin{lemma} \label{lemma:gammaq}
  Let $D_N$ be the Dirichlet kernel~\eqref{eq:7}.
  When $q>1$, the limit
  \begin{equation}
    \label{eq:9}
    \gamma(q) =
    \lim_{N \to \infty} \frac{\norm{D_N}_{L^q}}{N^{1-1/q}}
  \end{equation}
  exists, is finite, and its value is
  \begin{equation*}
    \gamma(q) =
    \tonde{\1\pi \int_{-\infty}^{+\infty}
      \abs{\frac{\sin x}{x}}^q \d x}^{1/q}.
  \end{equation*}
  Moreover, if $q>2$ then $\gamma(q) < 1$.
\end{lemma}  

\begin{remark} \label{rem:2m}
  When $q=2m$ is an even integer and $p=\infty$,
  it is easy to see that example~\ref{exa:dirichlet}
  provides a maximizer and hence
  \begin{equation*}
    C(\infty \to 2m) = \norm{D_N}_{L^{2m}}.
  \end{equation*}
  Indeed, this is an immediate consequence of
  the following monotonicity property.
  Let us assume that $\abs{a_n} \le b_n$
  for all $n=0,\dots,N-1$,
  and let $f = T_N(a)$ and $g = T_N(b)$,
  then, because of the positivity of the delta function, we have
  \begin{multline*}
    \norm{f}^{2m}_{L^{2m}} = 
    \1{2\pi} \int_{-\pi}^\pi
    \underbrace{f(t) \ovl{f(t)} \cdots
      f(t) \ovl{f(t)}}_{\text{$m$ times}} \d t = \\
    = \sum_{\substack{n_1, \dots, n_m \\
        \ovl{n}_1, \dots, \ovl{n}_m}}
    a_{n_1} \ovl{a_{\ovl{n}_1}} \cdots
    a_{n_m} \ovl{a_{\ovl{n}_m}}
    \cdot \1{2\pi} \int_{-\pi}^\pi
    \Eu^{\iC t \sum_j \tonde{n_j - \ovl{n}_j}} \d t = \\
    = \sum_{\substack{n_1, \dots, n_m \\
        \ovl{n}_1, \dots, \ovl{n}_m}}
    a_{n_1} \ovl{a_{\ovl{n}_1}}
    \cdots
    a_{n_m} \ovl{a_{\ovl{n}_m}}
    \cdot \delta\tonde{\sum_j \tonde{n_j - \ovl{n}_j}} \le \\
    \le \sum_{\substack{n_1, \dots, n_m \\
        \ovl{n}_1, \dots, \ovl{n}_m}}
    b_{n_1} b_{\ovl{n}_1}
    \cdots
    b_{n_m} b_{\ovl{n}_m}
    \cdot \delta\tonde{\sum_j \tonde{n_j - \ovl{n}_j}}
    = \norm{g}^{2m}_{L^{2m}}.
  \end{multline*} 
\end{remark}

\begin{remark}
  Using interpolation,
  from lemma~\ref{lemma:gammaq} and remark~\ref{rem:2m}
  it follows that we must have strict inequality,
  $C_N(p \to q) < N^{1-1/q-1/p}$,
  and also that
  \begin{equation*}
    \limsup_{N \to \infty}
    \frac{C_N(p \to q)}{N^{1 - 1/q - 1/p}} \le
    \gamma(q)^{1 - q'/p} < 1,
  \end{equation*}
  in the interior of region $\mB$
  where $2<q<\infty$ and $1/p < 1 - 1/q$.
\end{remark}

\subsection{Lower bounds for region $\mC$}

We now need an example for region $\mC$
which could give us a lower bound of the type
\begin{equation*}
  C_N(p \to q) \ges N^{1/2 - 1/p}.
\end{equation*}
A good candidate would be a choice of $a \in \C^N$
with $\abs{a_n} \sss 1$
for most of the $n$'s,
and such that $\abs{f(t)} \ges N^{1/2}$
for most of the $t$'s.

The idea is to set $a_n = \Eu^{\iC \varphi(n)}$,
for some real valued function $\varphi$,
and compare the sum
$\sum_{n = 0}^{N-1} \Eu^{\iC \varphi(n)} \Eu^{\iC n t}$
with the integral
$\int_0^N \Eu^{\iC \varphi(x)} \Eu^{\iC x t} \d x$,
with the help of the following lemma
taken from Zygmund's ``Trigonometric series'' \cite{Zyg1959}.
For the sake of completeness,
we present the interesting proof of the lemma 
in section~\ref{sec:proof-lemmata}.

\begin{lemma}[{\cite[chapter V, lemma 4.4]{Zyg1959}}] 
  \label{lemma:1}
  Let $\Phi$ be a smooth real valued function
  such that $\Phi'$ is monotone and
  \begin{equation*}
    \abs{\Phi'(t)} \le M < 2\pi.
  \end{equation*}
  Let $N$ be a positive integer,
  set $S = \sum_{n=1}^N \Eu^{\iC \Phi(n)}$
  and $I = \int_0^N \Eu^{\iC \Phi(x)} \d x$.
  Then $\abs{S-I} \le C_M$,
  where $C_M$ is a constant which depends only on $M$
  and does not depend on $\Phi$ or $N$.
\end{lemma}

We also need another lemma
whose proof is given in section~\ref{sec:proof-lemmata}.

\begin{lemma} \label{lemma:2}
  Let
  \begin{equation*}
    I(N,t) = \int_0^1 \Eu^{-\iC N (y-t)^2} \d y.
  \end{equation*}
  Then if $t \in ]0,1[$ we have
  \begin{equation*}
    \abs{I(N,t) - \sqrt{\frac{\pi}{N}} \Eu^{-\iC \pi/4}}
    \le \tonde{\1t + \1{1-t}} \1N.
  \end{equation*}
\end{lemma}

We are now ready to construct our example for region $\mC$.

\begin{example} 
  Let $\varphi(x) = -x^2/N$.
  Let us choose $a \in \C^N$ defined by 
  \begin{equation*}
    a_n = \Eu^{\iC\varphi(n)} = \Eu^{-\iC n^2/N}.
  \end{equation*}
  We have $\norm{a}_{\ell^p} = N^{1/p}$.
  We fix $t \in [-\pi, \pi]$
  and set $\Phi(x) = \varphi(x) + xt$.
  We have
  \begin{equation}
    \label{eq:10}
    f(t) = \sum_{n=0}^{N-1} \Eu^{\iC\varphi(n)} \Eu^{\iC nt} =
    \sum_{n=0}^{N-1} \Eu^{\iC\Phi(n)}
  \end{equation}
  The phase function $\Phi$ satisfies
  the hypotheses of lemma~\ref{lemma:1}:
  $\Phi'(x) = t - 2x/N$ is decreasing and
  for $0 \le x \le N$ we have
  \begin{equation*}
    \abs{\Phi'(x)} \le M = \pi + 2 < 2\pi.
  \end{equation*}
  It follows that the difference between the sum~\eqref{eq:10}
  and the integral
  \begin{equation*}
    g(t) = \int_0^N \Eu^{\iC\varphi(x)} \Eu^{\iC xt} \d x
  \end{equation*}
  is bounded by an absolute constant
  (independent of $N$ and $t$),
  \begin{equation}
    \label{eq:11}
    \abs{f(t) - g(t)} \le C.
  \end{equation}
  We have 
  \begin{equation*}
    g(t) = N \int_0^1 \Eu^{\iC N \tonde{yt - y^2}} \d y
    = N \Eu^{\iC Nt^2/4} \int_0^1 \Eu^{-\iC N(y - t/2)^2} \d y.
  \end{equation*}
  When $0<t<2$ we apply lemma~\ref{lemma:2} and obtain
  \begin{equation}
    \label{eq:12}
    \abs{g(t) - \sqrt\pi \Eu^{\iC(Nt^2 - \pi)/4} N^{1/2}} \le
    \frac{4}{t(2-t)}.
  \end{equation}
  From~\eqref{eq:11} and~\eqref{eq:12} we infer that
  \begin{equation*}
    \abs{f(t) - \sqrt\pi \Eu^{\iC(Nt^2 - \pi)/4} N^{1/2}} \le
    C + \frac{4}{t(2-t)},
  \end{equation*}
  for $0<t<2$.
  It follows that there exist
  two positive constants $0<\delta<1$ and $0<\eta<\sqrt\pi$,
  which do not depend on $N$,
  such that $\abs{f(t)} \ge \eta N^{1/2}$
  when $\delta < t < 2 - \delta$ and $N$ is sufficiently large.
  In particular we have
  \begin{equation*}
    \norm{f}_{L^q} \ge \tonde{2-2\delta}^{1/q} \eta N^{1/2},
  \end{equation*}
  which implies
  \begin{equation*}
    C_N(p \to q) \ge \frac{\norm{f}_{L^q}}{\norm{a}_{\ell^p}}
    \ge \eta \tonde{2-2\delta}^{1/q} N^{\12 - \1p}.
  \end{equation*}
  This shows that we can take
  $c_\mC = \min\graffe{\eta, \, \eta \sqrt{2 -2\delta}}$
  in theorem~\ref{theorem:1}.
\end{example}

\section{An integral operator with oscillating kernel} 
\label{sec:an-integral-operator}

We turn our attention to the linear integral operator
\begin{equation*}
  T_N: L^p([0,1];\C) \to L^q([0,1];\C)
\end{equation*}
defined by
\begin{equation*}
  T_Nf(t) = \int_0^1 \Eu^{\iC N / (1+t+s)} 
  \frac{f(s)}{(1+t+s)^\gamma} \d s,
\end{equation*}
for some fixed $\gamma \ge 0$.
Let us denote now by $C_N(p \to q)$
the best constant which can appear
in the estimate
$\norm{T_Nf}_{L^q} \le C \norm{f}_{L^p}$, \ie\
\begin{equation*}
  C_N(p \to q) = \norm{T_N}_{L^p \to L^q} =
  \sup_{f \in L^p \setminus 0}
  \frac{\norm{T_Nf}_{L^q}}{\norm{f}_{L^p}},
\end{equation*}
where this time
\begin{align*}
  \norm{f}_{L^p} &= \tonde{\int_0^1 \abs{f(t)}^p \d t}^{1/p},\;
  1 \le p < \infty, &
  \norm{f}_{L^\infty} &= \esssup_{t \in [0,1]} \abs{f(t)}.
\end{align*}
We ask the same question as before:
what can we say about the behaviour of $C_N$ as $N \to \infty$?

\begin{theorem} \label{theorem:2}
There exists positive absolute constants
$c_{\mA}$, $c_{\mB}$ and $c_{\mC}$
such that:
  \begin{align*}
    (p \to q) \in \mA &\implies
    c_{\mA}^{-1} (1+N)^{-1+1/p} \le
    C_N(p \to q) \le c_{\mA} (1+N)^{-1+1/p}; \\
    (p \to q) \in \mB &\implies
    c_{\mB}^{-1} (1+N)^{-1/q} \le
    C_N(p \to q) \le c_{\mB} (1+N)^{-1/q}; \\
    (p \to q) \in \mC &\implies
    c_{\mC}^{-1} (1+N)^{-1/2} \le
    C_N(p \to q) \le c_{\mC} (1+N)^{-1/2}.
  \end{align*}
\end{theorem}

In particular, the theorem says that
\begin{equation}
  \label{eq:13}
  \norm{T_N}_{L^p \to L^q} \les
  (1+N)^{-\min\graffe{1-1/p,\; 1/q ,\; 1/2}}.
\end{equation}

The following proof is similar
to the proof of theorem~\ref{theorem:1}.

\subsection{Upper bounds}

In order to prove the upper bounds in theorem~\ref{theorem:2},
it is enough to observe that we have
the $(1 \to \infty)$ \emph{dispersive} estimate
\begin{equation}
  \label{eq:14}
  \norm{T_Nf}_{L^\infty} \le \norm{f}_{L^1},
\end{equation}
and the $(2 \to 2)$ \emph{energy} estimate
\begin{equation}
  \label{eq:15}
  \norm{T_Nf}_{L^2} \les (1+N)^{-1/2} \norm{f}_{L^2}.
\end{equation}
Interpolation yields the $(p \to p')$ estimate
\begin{equation*}
  \norm{T_Nf}_{L^{p'}} \les (1+N)^{-1+1/p} \norm{f}_{L^p},
\end{equation*}
when $1 \le p \le 2$,
and H\"older's inequality does the rest.

Estimate~\eqref{eq:14} is trivial
and simply follows
by taking absolute values inside the integral.

To get estimate~\eqref{eq:15},
let $\chi \in C_0^\infty(\R)$ be a non-negative cut-off function
such that $\chi(t) = 1$ for $t \in [0,1]$
and $\chi(t) = 0$ for $t \notin [-1/2, 3/2]$;
then
\begin{equation*}
  \norm{T_Nf}_{L^2}^2 \le
  \int \chi(t) T_Nf(t) \ovl{T_Nf(t)} \d t = 
  \int_0^1 \int_0^1 K_N(s,\sigma)
  f(s) \ovl{f(\sigma)} \d s \d\sigma,
\end{equation*}
where the kernel $K_N$ is given by the oscillatory integral
\begin{equation*}
  K_N(s,\sigma) =
  \int \Eu^{\iC N \varphi(t;s,\sigma)}
  \chi(t;s,\sigma) \d t
\end{equation*}
with phase
\begin{equation*}
   \varphi(t;s,\sigma) = \1{1+t+s} - \1{1+t+\sigma} = 
   \frac{\sigma - s}{(1+t+s)(1+t+\sigma)},
\end{equation*}
and amplitude
\begin{equation*}
  \chi(t;s,\sigma) = \frac{\chi(t)}
  {(1+t+s)^\gamma (1+t+\sigma)^\gamma}.
\end{equation*}
To estimate $K_N$ we apply the principle of non-stationary phase
as illustrated by the following lemma.

\begin{lemma}[{\cite[chapter VIII]{Ste1993}}]
  \label{lemma:3}
  Let $K \in \N$.
  Let $\chi \in C_0^K(]a,b[)$
  and $\psi \in C^{K+1}(]a,b[)$.
  Suppose there exist positive constants $M, \delta$
  such that $\abs{\chi^{(k)}(t)} \le M$,
  $\abs{\psi^{(k+1)}(t)} \le M$
  for $k=0, \dots, K$ and $t \in ]a,b[$,
  and  $\abs{\psi'(t)} \ge \delta$.
  For $\lambda \in \R$, define
  \begin{equation*}
    I(\lambda) =
    \int_a^b \Eu^{\iC \lambda \psi(t)} \chi(t) \d t.
  \end{equation*}
  Then
  \begin{equation*}
    \abs{I(\lambda)} \le
    \frac{C(K, M, a, b)}
    {\tonde{1 + \delta \abs{\lambda}}^K}.
  \end{equation*}
\end{lemma}

When $s,\sigma \in [0,1]$, $t \in [-1/2, 3/2]$ and $k \in \N$, 
we have the following bounds
for the phase $\varphi(t;s,\sigma)$,
the amplitude $\chi(t;s,\sigma)$
and their derivatives:
\begin{align*}
   & C_k^{-1} \abs{s-\sigma} \le
   \abs{\de_t^k \varphi(t;s,\sigma)}
   \le C_k \abs{s-\sigma}, &
   & \abs{\de_t^k \chi(t;s,\sigma)} \le C_k.
\end{align*}
We apply lemma~\ref{lemma:3} with $K=2$,
$\lambda = N \abs{s-\sigma}$,
$\psi(t) = \varphi(t;s,\sigma)$,
so that
\begin{equation*}
  \abs{K_N(s,\sigma)} \les \1{\tonde{1 + N \abs{s-\sigma}}^2}.
\end{equation*}
Finally, we obtain~\eqref{eq:15} using Young's inequality,
\begin{equation*}
  \norm{T_Nf}_{L^2}^2 \les
  \int_0^1 \int_0^1
  \frac{\abs{f(s)} \cdot \abs{f(\sigma)}}
  {\tonde{1 + N \abs{s-\sigma}}^2}
  \d s \d\sigma \les \1{1+N} \norm{f}_{L^2}^2,
\end{equation*}
since we have $\norm{1/(1+N\abs{s})^2}_{L^1(\R)} \les 1/(1+N)$.

\subsection{Lower bounds}

The lower bounds in theorem~\ref{theorem:2},
are a consequence of the following examples
which are inspired by the examples
of section~\ref{sec:lower-bounds}.

\begin{example}
  Let $\eta$ be a small positive constant.
  Let $f(s) = \chi\tonde{0 \le s \le \eta/N}$
  be the characteristic function of the interval $[0,\eta/N]$.
  When $t \in [0,1]$ and $s \in [0,\eta/N]$ we have
  \begin{equation*}
    \frac{N}{1+t+s} = \frac{N}{1+t} + O(\eta),
  \end{equation*}
  so that
  \begin{equation*}
    \Eu^{\iC N/(1+t+s)} =
    \Eu^{\iC N/(1+t)} + O(\eta),
  \end{equation*}
  and
  \begin{equation*}
    T_Nf(t) =
    \Eu^{\iC N/(1+t)}
    \int_0^{\eta/N} \frac{\d s}{(1+t+s)^\gamma}
    + O\tonde{\frac{\eta^2}{N}} =
    \Eu^{\iC N/(1+t)} \cdot
    \1{(1+t)^\gamma} \cdot \frac{\eta}{N}
    + O\tonde{\frac{\eta^2}{N}}.
  \end{equation*}
  It follows that,
  if $\eta$ is sufficiently small,
  $\abs{T_Nf(t)} \ges N^{-1}$ 
  uniformly for $t\in [0,1]$.
  Hence, 
  \begin{equation*}
    C_N(p \to q) \ge \frac{\norm{T_N f}_{L^q}}{\norm{f}_{L^p}}
    \ges \frac{N^{-1} \cdot 1^{1/q}}{N^{-1/p}} = N^{-1+1/p}.
  \end{equation*}
\end{example}

\begin{example}
  Let $f(s) = e^{-\iC N/(1+s)}$ for $s \in [0,1]$.
  Let $\eta$ be a small positive constant.
  When $s \in [0,1]$ and $t \in [0,\eta/N]$ we have
  \begin{equation*}
    \frac{N}{1+t+s} - \frac{N}{1+s} = O(\eta),
  \end{equation*}
  so that
  \begin{equation*}
    \Eu^{\iC N/(1+t+s)} \Eu^{-\iC N/(1+s)}
    = \tonde{1 + O(\eta)},
  \end{equation*}
  and
  \begin{equation*}
    T_Nf(t) = \int_0^1 \frac{\d s}{(1+t+s)^\gamma} + O(\eta).
  \end{equation*}
  It follows that,
  if $\eta$ is sufficiently small,
  $\abs{T_Nf(t)} \ges 1$ 
  uniformly for $t\in [0,\eta/N]$.
  Hence, 
  \begin{equation*}
    C_N(p \to q) \ge \frac{\norm{T_N f}_{L^q}}{\norm{f}_{L^p}}
    \ges \frac{1 \cdot N^{-1/q}}{1^{1/p}} = N^{-1/q}.
  \end{equation*}
\end{example}

\begin{example}
  Let $f(s) = \Eu^{\iC N s^2}$.
  We have $T_Nf(t) = \int_0^1 \Eu^{\iC N \varphi(s;t)} \d s$,
  where the phase function is given by
  \begin{equation*}
    \varphi(s;t) = s^2 + \1{1+t+s}.
  \end{equation*}
  For each $t \in [0,1]$,
  the phase possesses exactly one non degenerate critical point,
  since its first derivative,
  \begin{equation*}
    \de_s \varphi(s;t) = 2s - \1{(1+t+s)^2},
  \end{equation*}
  vanishes in correspondence with the point $s_*$
  defined implicitely  by the equation
  \begin{equation*}
    s_* = \1{2(1+t+s_*)^2} \in \quadre{\1{18}, \12},
  \end{equation*}
  and its second derivative never vanishes,
  \begin{equation*}
    \de_s^2 \varphi(s;t) = 2 + \frac{2}{(1+t+s)^3} \ge 2.
  \end{equation*}
  From an application of
  the principle of stationary phase,
  lemma~\ref{lemma:statphase} below,
  it follows that $\abs{T_Nf(t)} \sss N^{-1/2}$
  uniformly for $t \in [0,1]$.
  Hence,
  \begin{equation*}
    C_N(p \to q) \ge \frac{\norm{T_N f}_{L^q}}{\norm{f}_{L^p}}
    \ges \frac{N^{-1/2} \cdot 1^{1/q}}{1^{1/p}} = N^{-1/2}.
  \end{equation*}
\end{example}

\begin{lemma}[{\cite[chapter VIII]{Ste1993}}]
  \label{lemma:statphase}
  Let us consider the oscillatory integral
  \begin{equation}
    \label{eq:16}
    I(N) = \int_a^b \Eu^{\iC N \varphi(s)} \chi(s) \d s,
  \end{equation}
  with phase function $\varphi \in C^5([a,b])$
  and amplitude $\chi \in C^2([a,b])$.
  We assume that $\varphi''(s) \ge 1$ for all $s \in [a,b]$
  and that $\varphi'(s_*) = 0$
  for a point $s_* \in [a+\delta,b-\delta]$,
  with $\delta>0$.
  Then
  \begin{equation*}
    I(N) = \frac{J_* \Eu^{\iC N \varphi(s_*)}}{\sqrt{N}}
    + O\tonde{\1N},
  \end{equation*}
  where
  \begin{equation*}
    J_* = \Eu^{\iC \pi/4}
    \chi(s_*) \sqrt{\frac{2\pi}{\varphi''(s_*)}},
  \end{equation*}
  and the implicit constant in the $O$-symbol
  depends only on the quantities $b-a$, 
  $\varphi(a)$, $\varphi(b)$, $\delta$
  and on uniform bounds for
  $\abs{\chi^{(k)}}$ and $\abs{\varphi^{(k+3)}}$,
  with $k=0,1,2$.
\end{lemma}

\section{Local estimates for
  inhomogeneous Schr\"odinger equations}
\label{sec:local-estim-inhom}

Now we come to the problem which motivated the above study,
namely the problem of determining
the optimal range of exponents
for local inhomogeneous Strichartz estimates.

Let $n \ge 3$.
Let $u(t,x)$ be the solution of
the inhomogeneous Schr\"odinger equation
\begin{align*}
  i\de_t u - \Delta u &= F(t,x), &
  &(t,x) \in \R \times \R^n,
\end{align*}
with zero initial data $u(0,x) = 0$.
We assume that the support of $F$ is contained
in the region where $0 \le t \le 1$
and we look at the solution $u(t,x)$ 
in the region where $2 \le t \le 3$.
Using the fundamental solution of the Schr\"odinger equation
we can write an explicit formula for $u$ in terms of $F$:
\begin{equation}
  \label{eq:17}
  u(t,x) = (4\pi)^{-n/2} \int_0^1 \int_{\R^n}
  \frac{\Eu^{-\iC |x-y|^2/(4(t-s))}}
  {(t-s)^{n/2}} F(s,y) \d y \d s.
\end{equation}

\begin{figure}
  \centering
  \begin{picture}(0,0)%
    \includegraphics{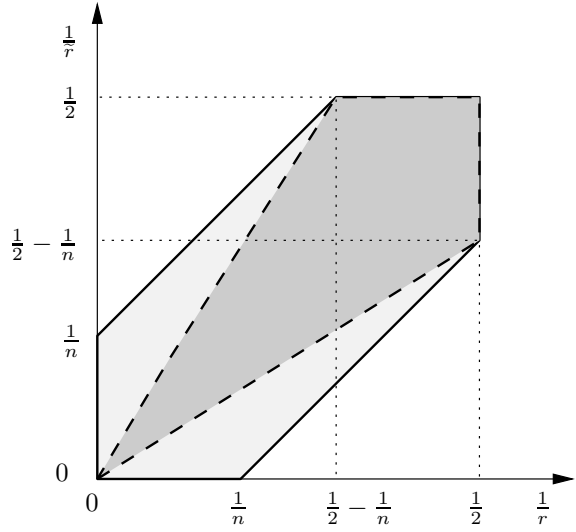}%
  \end{picture}%
  \setlength{\unitlength}{3947sp}%
  \begin{picture}(3567,3286)(46,-2435)
    \put(2926,-2386){\makebox(0,0)[lb]{\smash{$\12$}}}
    \put(2026,-2386){\makebox(0,0)[lb]{\smash{$\12-\1n$}}}
    \put(1426,-2386){\makebox(0,0)[lb]{\smash{$\1n$}}}
    \put(376,-1336){\makebox(0,0)[lb]{\smash{$\1n$}}}
    \put(376,164){\makebox(0,0)[lb]{\smash{$\12$}}}
    \put(338,-2176){\makebox(0,0)[lb]{\smash{$0$}}}
    \put(526,-2363){\makebox(0,0)[lb]{\smash{$0$}}}
    \put( 46,-729){\makebox(0,0)[lb]{\smash{$\12-\1n$}}}
    \put(3345,-2386){\makebox(0,0)[lb]{\smash{$\1r$}}}
    \put(369,547){\makebox(0,0)[lb]{\smash{$\1\rt$}}}
  \end{picture}
  \caption{The region of admissible exponents
    for estimates~\eqref{eq:18} and~\eqref{eq:21}
    in the $(1/r, 1/\rt)$ plane.}
  \label{fig:rrt}
\end{figure}

Local Strichartz estimates of the type
\begin{equation}
  \label{eq:18}
  \norm{u}_{L^q([2,3]; L^r(\R^n))} \les
  \norm{F}_{L^{\qt'}([0,1]; L^{\rt'}(\R^n))}
\end{equation}
are known to hold when
the pairs of exponents $(q,r)$ and $(\qt, \rt)$
satisfy the conditions
\begin{align*}
  &\12 - \1n \le \1r \le \12, &
  &\1q \ge \nd \tonde{\12 - \1r}, \\
  &\12 - \1n \le \1\rt \le \12, &
  &\1\qt \ge \nd \tonde{\12 - \1\rt},  
\end{align*}
(see~\cite{Kat1994} and~\cite{KeeTao1998}
for details and references).
The norms which appear in~\eqref{eq:18}
have the integration with respect to space variables
computed before doing the integration with respect to time.
Interpolation with the easy dispersive estimate
\begin{equation*}
  \norm{u}_{L^\infty([2,3]\times\R^n)} \les
  \norm{F}_{L^1([0,1]\times\R^n)},
\end{equation*}
proves that when exponents $r$ and $\rt$ satisfy the conditions
\begin{align}
  \label{eq:19}
  & \frac{n-2}{n} \cdot \1r \le \1\rt \le \12, &
  & \frac{n-2}{n} \cdot \1\rt \le \1r \le \12,  
\end{align}
then there exist some exponents $q, \qt \in [1, \infty]$
for which estimate~\eqref{eq:18} holds~(\cite{Fos2003pre}).
The dark shaded area in figure~\ref{fig:rrt}
shows the region in the $(r,\rt)$ plane
where~\eqref{eq:19} is satisfied.

Here,
we use the result of section~\ref{sec:an-integral-operator}
to obtain estimates similar to~\eqref{eq:18}
but with norms which have the integration with respect to time 
computed before doing the integration with respect 
to the space variables.

\begin{theorem} \label{theorem:3}
  Let us suppose that exponents $r$ and $\rt$
  satisfy the conditions
  \begin{align}
    \label{eq:20}
    & 0 \le \1r \le \12, &
    & 0 \le \1\rt \le \12, &
    & \abs{\1r - \1\rt} < \1n,
  \end{align}
  then there exist some exponents $q, \qt \in [1, \infty]$
  for which we have the estimate
  \begin{equation}
    \label{eq:21}
    \norm{u}_{L^r(\R^n; L^q([2,3]))} \les
    \norm{F}_{L^{\rt'}(\R^n; L^{\qt'}([0,1]))}.
  \end{equation}
\end{theorem}

The light shaded area in figure~\ref{fig:rrt}
shows the region in the $(r,\rt)$ plane
where~\eqref{eq:20} is satisfied.

\begin{remark} \label{rem:weaker}
  Estimate~\eqref{eq:21} is weaker than~\eqref{eq:18}
  in the sense that when~\eqref{eq:18} holds
  for the pairs of exponents $(q,r)$ and $(\qt, \rt)$
  then~\eqref{eq:21} holds if we replace
  $q$ with $\min\!\graffe{q,r}$ and
  $\qt$ with $\min\!\graffe{\qt, \rt}$.
  Indeed, if $Q = \min\!\graffe{q,r}$, then we have
  \begin{equation*}
    \norm{G}_{L^r(\R^n;L^Q(I))} \le
    \norm{G}_{L^Q(I;L^r(\R^n))} \le
    \norm{G}_{L^q(I;L^r(\R^n))},
  \end{equation*}
  and
  \begin{equation*}
    \norm{G}_{L^{q'}(I;L^{r'}(\R^n))} \le
    \norm{G}_{L^{Q'}(I;L^{r'}(\R^n))} \le
    \norm{G}_{L^{r'}(\R^n;L^{Q'}(I))},
  \end{equation*}
  for any function $G(t,x)$ and time interval $I$ of length $1$.
\end{remark}

\begin{proof}[Proof of theorem~\ref{theorem:3}]
  Using the change of variables
  \begin{align*}
   & t \in [2,3] \changeto \tau=t-2 \in [0,1], \\
   & s \in [0,1] \changeto \sigma=1-s \in [0,1],
  \end{align*}
  we have $t-s = 1+\tau+\sigma$ and we see that 
  \begin{equation*}
    u(t,x) = (4\pi)^{-n/2} \int_{\R^n} \int_0^1
    \Eu^{-\iC (|x-y|^2/4)/(1+\tau+\sigma)}
    \frac{F(1-\sigma,y)}{(1+\tau+\sigma)^{n/2}} \d \sigma \d y.
  \end{equation*}
  The inner integral is exacly the operator $T_N$
  described in section~\ref{sec:an-integral-operator},
  with $N = |x-y|^2/4$ and $\gamma = n/2$, 
  acting on the function $\sigma \mapsto F(1-\sigma, y)$.
  We can apply theorem~\ref{theorem:2} and
  use~\eqref{eq:13} to obtain
  \begin{equation}
    \label{eq:22}
    \norm{u(\cdot,x)}_{L^q([2,3])} \les
    \int_{\R^n} \frac{\norm{F(\cdot, y)}_{L^{\qt'}([0,1])}}
    {\tonde{1+\abs{x-y}}^\alpha} \d t,
  \end{equation}
  where
  \begin{equation*}
    \alpha = \min\graffe{\frac{2}{\qt},\; \frac{2}{q},\; 1}.
  \end{equation*}
  Estimate~\eqref{eq:21} then follows from~\eqref{eq:22}
  and Young's inequality
  when $1/r + 1/\rt < \alpha/n$,
  or the Hardy-Littlewood-Sobolev inequality
  (\cite[chapter VIII, section 4.2]{Ste1993})
  when $1/r + 1/\rt = \alpha/n$ and $r,\rt < \infty$.
  This means that we proved estimate~\eqref{eq:21} when
  \begin{equation}
    \label{eq:23}
    \1r + \1{\rt} \le \1n, \qquad
    \1q \ge \nd \tonde{\1r + \1{\rt}}, \qquad
    \1{\qt} \ge \nd \tonde{\1r + \1{\rt}},
  \end{equation}
  with strict inequalities if $r=\infty$ or $\rt=\infty$.
  On the other hand, 
  it follows from remark~\ref{rem:weaker}
  that estimate~\eqref{eq:21} holds for some $q$ and $\qt$
  when $r$ and $\rt$ satisfy conditions~\eqref{eq:19}.
  Using interpolation we obtain estimate~\eqref{eq:21}
  for some $q$ and $\qt$ whenever $r$ and $\rt$ 
  are in the convex hull of the two regions
  in the $(1/r, 1/\rt)$ plane described
  by~\eqref{eq:19} and~\eqref{eq:23}.
  This convex hull is precisely the region
  described by~\eqref{eq:20}
  plus the two endpoints $r=2$, $\rt = 2n/(n-2)$
  and $r=2n/(n-2)$, $\rt=2$.
\end{proof}

\begin{remark}
  The range~\eqref{eq:20} for exponents $r$ and $\rt$
  which appear in estimate~\eqref{eq:21}
  is almost optimal.
  In fact example~\ref{exa:rrtopt}
  shows that if estimate~\eqref{eq:21} holds,
  then we must necessarily have
  $1/r - 1/\rt \le 1/n$.

  The dual of the operator defined by~\eqref{eq:17},
  which takes $F \in L^{\qt'}([0,1]; L^{\rt'}(\R^n))$
  into $u \in L^q([2,3]; L^r(\R^n))$,
  is an operator with the same structure.
  Hence, by duality we must also have
  $1/\rt - 1/r \le 1/n$.
\end{remark}

\begin{example} \label{exa:rrtopt}
  Let $R \gg 1$ and $0< \eta <1$.
  We choose
  \begin{equation*}
    F(s,y) = \Eu^{\iC 2 R^2 s^2} \chi\tonde{Ry/\eta},
  \end{equation*}
  where $\chi$ is the characteristic function
  of the unit ball in $\R^n$ centered at the origin.
  We have 
  $\norm{F}_{L^{\rt'}(\R^n; L^{\qt'}([0,1]))} \sss (\eta/R)^n$.
  We can write the solution $u$ as
  \begin{equation}
    \label{eq:38}
    u(t,x) = \int_{|y|<\eta/R} \int_0^1
    \frac{\Eu^{\iC \tonde{2 R^2 s^2 - |x-y|^2/(4(t-s))}}}
    {(t-s)^{n/2}} \d s \d y
    = \int_{|y|<\eta/R} I\tonde{t, \frac{x-y}{2R}, R} \d y,
  \end{equation}
  where $I$ is the oscillatory integral
  \begin{equation*}
    I(t,z,R) =
    \int_0^1 \Eu^{\iC R^2 \varphi(s;t,z)} \psi(s;t) \d s,
  \end{equation*}
  with phase $\varphi(s;t,z) = 2s^2 - |z|^2/(t-s)$
  and amplitude $\psi(s;t) = 1/(t-s)^{n/2}$.
  The first and second derivatives of the phase are
  \begin{align*}
    \de_s \varphi &= 4s - \frac{|z|^2}{(t-s)^2}, &
    \de_s^2 \varphi &= 4 - \frac{2|z|^2}{(t-s)^3}.
  \end{align*}
  When $t \in [2,3]$ and $1/2 \le |z| \le 1$,
  derivatives with respect to $s$ of all orders
  for $\varphi$ and $\psi$
  are uniformly bounded by absolute constants
  and the phase $\varphi$
  has exactly one non degenerate critical point
  $s_* = s_*(t,z)$ in $[0,1]$,
  \begin{align*}
    s_* &= \frac{|z|^2}{4(t-s_*)^2} \in \quadre{\1{48}, \14}, &
    \de_s \varphi(s_*, t, z) &= 0, &
    \de_s^2 \varphi(s_*, t, z) &= 4 - \frac{2s_*}{t-s_*}
    \in \quadre{\frac23,12}.
  \end{align*}
  We are in the conditions to apply lemma~\ref{lemma:statphase}
  and obtain
  \begin{equation*}
    I(t,z,R) =
    \frac{J_*(t,z) \Eu^{\iC R^2 \varphi_*(t,z)}}{R}
    + O\tonde{\1{R^2}},
  \end{equation*}
  where
  \begin{align*}
    J_*(t,z) &= \Eu^{i\pi/4} \psi(s_*(t,z);t)
    \sqrt{\frac{2\pi}{\de_s^2 \varphi(s_*(t,z);t,z)}}, &
    \varphi_*(t,z) &= \varphi(s_*(t,z);t,z).
  \end{align*}
  By the above computations, when $t \in [2,3]$ and $|z|<1$
  we have
  \begin{equation}
    \label{eq:39}
    \abs{J_*(t,z)} \sss 1;
  \end{equation}
  moreover, 
  \begin{equation*}
    \nabla_z \varphi_*(t,z) = \nabla_z \varphi(s_*;t,z) =
    - \frac{2z}{t-s_*} = O(1),
  \end{equation*}
  so that we have
  $\varphi_*(t,z) - \varphi(t,z_0) = O(z-z_0)$.
  In particular this shows that the oscillatory factor
  \begin{equation*}
     \Eu^{\iC R^2 \varphi_*\tonde{t, \frac{x-y}{2R}}}
     = \Eu^{\iC R^2 \varphi_*\tonde{t, \frac{x}{2R}} + O(\eta)}
     = \Eu^{\iC R^2 \varphi_*\tonde{t, \frac{x}{2R}}}
     \tonde{1+O(\eta)},
  \end{equation*}
  does not oscillates too much 
  when $|y| \le \eta/R$, $R \le |x-y| \le 2R$
  and $\eta$ is sufficiently small.
  It follows that
  \begin{equation*}
    I\tonde{t, \frac{x-y}{2R}, R} = 
    \frac{\Eu^{\iC R^2 \varphi_*(t,x/(2R))}}{R}
    J_*\tonde{t, \frac{x-y}{2R}} \tonde{1 + O(\eta)}
    + O\tonde{\1{R^2}},
  \end{equation*}
  which inserted in~\eqref{eq:38} and using~\eqref{eq:39}
  proves that
  \begin{equation*}
    \abs{u(t,x)} \ges \frac{\eta^n}{R^{n+1}},
  \end{equation*}
  on the region $2\le t\le3$, $R+\eta/R<|x|<2R-\eta/R$.
  Thus, the ratio
  \begin{equation*}
    \frac{\norm{u}_{L^r(\R^n; L^q([2,3]))}}
    {\norm{F}_{L^{\rt'}(\R^n; L^{\qt'}([0,1]))}} \ges
    \frac{R^{-n-1} \cdot R^{n/r}}{R^{-n/\rt'}} = 
    R^{n(1/r - 1/\rt) - 1},
  \end{equation*}
  cannot be bounded as $R \to \infty$ unless
  the necessary condition
  \begin{equation*}
    \1r - \1\rt \le \1n
  \end{equation*}
  is satisfied.
\end{example}

\section{Proof of the lemmata} \label{sec:proof-lemmata}

\begin{proof}[Proof of lemma~\ref{lemma:gammaq}]
  We want to compute the leading term 
  in the asymptotic espansion of the $L^q$
  norm of the Dirichlet kernel $D_N$ as $N \to \infty$.
  By a simple change of variable we have
  \begin{equation}
    \label{eq:25}
    \norm{D_N}_{L^q}^q = \1{N \pi} \int_{-N\pi/2}^{N \pi/2}
    \abs{\frac{\sin x}{\sin(x/N)}}^q \d x.
  \end{equation}
  When $\abs{\alpha}<\pi/2$ we have
  \begin{equation*}
    \abs{\1{\sin\alpha}}^q = \abs{\1\alpha}^q
    + O\tonde{\abs{\alpha}^{2-q}}.
  \end{equation*}
  It follows that
  \begin{multline}
    \label{eq:26}
    \int_{-N\pi/2}^{N \pi/2}
    \abs{\frac{\sin x}{\sin(x/N)}}^q \d x =
    N^q \int_{-N\pi/2}^{N \pi/2} 
    \abs{\frac{\sin x}{x}}^q \d x + {}\\
    + O\tonde{N^{q-2} \int_{-N\pi/2}^{N\pi/2}
      \abs{\sin x}^q \abs{x}^{2-q} \d x}.
  \end{multline}
  When $q>1$,
  we have $\int_N^\infty x^{-q} \d x = o(1)$,
  and $\int_1^N x^{2-q} \d x = o(N^2)$ as $N \to \infty$,
  hence
  \begin{equation}
    \label{eq:27}
    \1\pi
    \int_{-N\pi/2}^{N \pi/2} \abs{\frac{\sin x}{x}}^q \d x =
    \1\pi
    \int_{-\infty}^{+\infty} \abs{\frac{\sin x}{x}}^q \d x +
    o(1)
    = \gamma(q)^q + o(1).
  \end{equation}
  and
  \begin{equation}
    \label{eq:28}
    \int_{-N\pi/2}^{N\pi/2} \abs{\sin x}^q \abs{x}^{2-q} \d x =
    o(N^2).
  \end{equation}
  When we insert~\eqref{eq:26},~\eqref{eq:27} and~\eqref{eq:28}
  into~\eqref{eq:25}, we obtain
  \begin{equation*}
    \norm{D_N}_{L^q}^q = N^{q-1} \gamma(q)^q + o\tonde{N^{q-1}},
  \end{equation*}
  which proves~\eqref{eq:9}.
  The quantity $\gamma(q)^q$ is decreasing, 
  indeed
  \begin{equation*}
    \der{}{q} \quadre{\gamma(q)^q} =
    \1\pi \int_{-\infty}^{+\infty}
    \abs{\frac{\sin x}{x}}^q \log\abs{\frac{\sin x}{x}} \d x
    < 0.
  \end{equation*}
  Thus, $\gamma(q)^q < 1$ when $q>2$,
  since we know from remark~\ref{rem:energy}
  that $\gamma(2) = 1$.
\end{proof}

\begin{proof}[Proof of lemma~\ref{lemma:1}]
  Integration by parts gives
  \begin{equation*}
    I = N \Eu^{\iC \Phi(N)}
    - \int_0^N x \de_x \Eu^{\iC \Phi(x)} \d x;
  \end{equation*}
  summation by parts \cite[chapter I, section 2]{Zyg1959} gives
  \begin{equation*}
    S = N \Eu^{\iC \Phi(N)} -
    \sum_{n=0}^{N-1} n \tonde{\Eu^{\iC \Phi(n+1)}
      - \Eu^{\iC \Phi(n)}}.
  \end{equation*}
  We write $\Eu^{\iC\Phi(n+1)} - \Eu^{\iC \Phi(n)} =
  \int_n^{n+1} \de_x \Eu^{\iC\Phi(x)} \d x$ and obtain
  \begin{equation*}
    S-I = \int_0^N \tonde{x - [x]} \de_x \Eu^{\iC\Phi(x)} \d x,
  \end{equation*}
  where $[x]$ denotes the largest integer
  smaller or equal to $x$.
  The function $x-[x]$ is
  a piecewise-continuous periodic function with period $1$
  and mean value $1/2$.
  When $x \notin \Z$, it coincides almost everywhere
  with its Fourier series given by
  \begin{equation*}
    x - [x] =
    \12 + \sum_{k \ne 0} \frac{\Eu^{2\pi \iC k x}}{2\pi \iC k}.
  \end{equation*}
  Thus, we have
  \begin{equation}
    \label{eq:29}
    S-I = \12 \tonde{\Eu^{\iC \Phi(N)} - \Eu^{\iC \Phi(0)}} + 
    \1{2\pi \iC} \sum_{k \ne 0} \frac{C(N,k)}{k},
  \end{equation}
  where
  \begin{multline}
    \label{eq:30}
    C(N,k) = \int_0^N \Eu^{2\pi \iC kx}
    \de_x \Eu^{\iC \Phi(x)} \d x
    = \int_0^N \iC\Phi'(x)
    \Eu^{\iC\tonde{\Phi(x)+2\pi kx}} \d x = \\
    = \int_0^N \frac{\Phi'(x)}{\Phi'(x)+2\pi k}
    \de_x \Eu^{\iC\tonde{\Phi(x)+2\pi kx}} \d x.
  \end{multline}
  Let $\Psi(x) = \Phi'(x) / (\Phi'(x)+2\pi k)$.
  We observe that the denominator does not vanish,
  since $k \ne 0$ and $|\Phi(x)|<2\pi$.
  The function $y \mapsto \Gamma(y) = y/(y+2\pi k)$ 
  is smooth and monotone on the interval $[-M, M]$ and 
  \begin{equation}
    \label{eq:31}
    \sup_{y \in [-M,M]} \abs{\Gamma(y)}
    \le \frac{M}{2\pi\abs{k} - M}
    \le \frac{M}{2\pi - M} \cdot \1{\abs{k}}.
  \end{equation}
  It follows that $\Psi = \Gamma \compose \Phi'$
  is also smooth and monotone,
  being the composition of two smooth and monotone functions.
  In particular $\de_x \Psi$ does not change sign and we have
  \begin{equation}
    \label{eq:32}
    \int_0^N \abs{\de_x \Psi(x)} \d x =
    \abs{\int_0^N \de_x \Psi(x) \d x} = 
    \abs{\Psi(N) - \Psi(0)}.
  \end{equation}
  Integrating~\eqref{eq:30} by parts we obtain
  \begin{equation*}
    C(N,k) =
    \quadre{\Psi(x) \Eu^{\iC\tonde{\Phi(x)+2\pi kx}}}_0^N
    - \int_0^N \tonde{\de_x \Psi(x)}
    \Eu^{\iC\tonde{\Phi(x)+2\pi kx}} \d x.
  \end{equation*}
  We now take absolute values
  and use~\eqref{eq:32} and~\eqref{eq:31},
  \begin{equation*}
    \abs{C(N,k)} \le
    \abs{\Psi(N)} + \abs{\Psi(0)} + \abs{\Psi(N) - \Psi(0)}
    \le \frac{4M}{2\pi-M} \cdot \1{\abs{k}}.
  \end{equation*}
  Substituting this in~\eqref{eq:29} we obtain
  \begin{equation*}
    \abs{S-I} \le
    1 + \1{2\pi}\cdot \frac{4M}{2\pi-M} \sum_{k \ne 0} \1{k^2}
    = C_M.
  \end{equation*}
\end{proof}

\begin{proof}[Proof of lemma~\ref{lemma:2}]
  Let $J(N, s) = \int_0^\infty \Eu^{-\iC N (z + s)^2} \d z$.
  If $s>0$, we write
  \begin{equation}
    \Eu^{-\iC N(z+s)^2} =
    - \1{2\iC N(z+s)} \de_z \Eu^{-\iC N(z+s)^2}
  \end{equation}
  and integrate by parts,
  \begin{equation}
    J(N,s) =
    \int_0^\infty
    \frac{\de_z \Eu^{-\iC N(z+s)^2}}{2iN(z+s)} \d z
    = \frac{\Eu^{-\iC Ns^2}}{2\iC Ns} -
    \int_0^\infty
    \frac{\Eu^{-\iC N(z+s)^2}}{2\iC N(z+s)^2} \d z.
  \end{equation}
  It follows that
  \begin{equation}
    \label{eq:33}
    \abs{J(N,s)} \le
    \1{2Ns} + \int_0^\infty \frac{\d z}{2N(z+s)^2} = \1{Ns}.
  \end{equation}
  We have
  \begin{equation*}
    \int_{-\infty}^{+\infty} \Eu^{-\iC N (y-t)^2} \d y
    - I(N,t) =
    J(N, t) + J(N, 1-t).
  \end{equation*}
  If $0<t<1$ the result follows from~\eqref{eq:33}
  and a rescaling of the Fresnel integral
  $\int_{-\infty}^{+\infty} \Eu^{-\iC x^2} \d x
  = \sqrt{\pi} \Eu^{-\iC \pi/4}$.
\end{proof}

\begin{proof}[Proof of lemma~\ref{lemma:3}]
  The method is standard.
  We begin by noticing that 
  \begin{equation*}
    \Eu^{\iC \lambda \psi(t)} = \1{\iC \lambda \psi'(t)} 
    \de_t \Eu^{\iC \lambda \psi(t)},
  \end{equation*}
  thus $\Eu^{\iC \lambda \psi}$ is an eigenfunction
  for the differential operator
  \begin{equation*}
    L f (t) = \1{\iC \psi'(t)} \de_t f(t),
  \end{equation*}
  relative to the eigenvalue $\lambda$,
  and in particular
  \begin{equation*}
    \Eu^{\iC \lambda \psi} =
    \lambda^{-1} L\tonde{\Eu^{\iC \lambda \psi}} = \dots =
    \lambda^{-K} L^K\tonde{\Eu^{\iC \lambda \psi}}.
  \end{equation*}
  Let $L^*$ be the transpose of $L$,
  \begin{equation*}
    L^* f (t) = - \de_t \tonde{\frac{f(t)}{\iC \psi'(t)}}.
  \end{equation*}
  Then, integrating by parts $k$ times we obtain
  \begin{equation*}
    \int \Eu^{\iC \lambda \psi} \chi \d t = \lambda^{-k}
    \int L^k\tonde{\Eu^{\iC \lambda \psi}} \chi \d t =
    \int \Eu^{\iC \lambda \psi} (L^*)^k(\chi) \d t,
  \end{equation*}
  thus, for $k=0, \dots, K$ we have
  \begin{equation*}
    \abs{I(\lambda)} \les \lambda^{-k}
    \int \abs{(L^*)^k(\chi)} \d t,
  \end{equation*}
  The quantity $(L^*)^k(\chi)$ is a linear combination of terms
  which are the ratio of products of
  $\chi$, its derivatives and derivatives of $\psi'$,
  with powers of $\psi'$; hence, its integral can be bounded
  by a constant which depends on $k$, $M$,
  and the measure of the support of $\chi$.
\end{proof}

\begin{proof}[Proof of lemma~\ref{lemma:statphase}]
  As a first step we prove the lemma in the special case
  of quadratic phase $\varphi(s) = k (s-s_*)^2$,
  for some $k \ge 1/2$.
  We write
  \begin{equation*}
    \chi(s) = \chi(s_*) + (s-s_*) \Psi(s),
  \end{equation*}
  so that
  \begin{equation}
    \label{eq:34}
    I(N) = \chi(s_*) \int_a^b \Eu^{\iC N k(s-s_*)^2} \d s +
    \1{2 \iC N k}
    \int_a^b \tonde{\de_s \Eu^{\iC N k (s-s_*)^2}} \Psi(s) \d s.
  \end{equation}
  As in the proof of lemma~\ref{lemma:2},
  let $J(N,t) = \int_0^\infty \Eu^{\iC N (t+z)^2} \d z$.
  Using estimate~\eqref{eq:33} we obtain
  \begin{multline*}
    \int_a^b \Eu^{\iC N k(s-s_*)^2} \d s =
    \int_{-\infty}^{+\infty} \Eu^{\iC N k s^2} \d s
    - J(Nk,s_*-a) - J(Nk,b-s_*) = \\
    = \sqrt{\frac{\pi}{Nk}} \Eu^{\iC \pi/4}
    + O\tonde{\1{N k \delta}}.
  \end{multline*}
  Integration by parts in the second integral
  in the right hand side of~\eqref{eq:34} gives
  \begin{equation*}
    \abs{\int_a^b \tonde{\de_s \Eu^{\iC Nk(s-s_*)^2}}
      \Psi(s) \d s}
    \le \abs{\Psi(b)} + \abs{\Psi(b)}
    + \int_a^b \abs{\Psi'(s)} \d s.
  \end{equation*}
  From the explicit formulas
  \begin{gather*}
    \Psi(s) = \frac{\chi(s) - \chi(s_*)}{s - s_*} =
    \1{s-s_*} \int_{s_*}^s \chi'(\sigma) \d\sigma, \\
    \Psi'(s) = \frac{\chi'(s) - \Psi(s)}{s - s_*} =
    \1{(s-s_*)^2} \int_{s_*}^s \int_\sigma^s
    \chi''(\rho) \d\rho \d\sigma,    
  \end{gather*}
  and the hypotheses on $\chi$, it follows that
  \begin{align*}
    \max \abs{\Psi} &\le \max \abs{\chi'} \le 1, &
    \max \abs{\Psi'} &\le \max \abs{\chi''} \le 1.   
  \end{align*}
  Hence,
  \begin{equation}
    \label{eq:35}
    I(N) = \chi(s_*) \sqrt{\frac{\pi}{Nk}} \Eu^{\iC \pi/4}
    + O\tonde{\frac{1+|b-a|}{N \delta}},
  \end{equation}
  which concludes the proof in this special case.

  In the second step, 
  we prove the lemma for a general phase function $\varphi$
  but we assume the amplitude $\chi$ to have compact support 
  in $]a,b[$.
  Let
  \begin{equation*}
    \Phi(s;\lambda) = 
    (1-\lambda) \tonde{\varphi(s_*)
      + \frac{\varphi''(s_*)}{2} (s-s_*)^2}
    + \lambda \varphi(s),
  \end{equation*}
  and $I(N;\lambda) =
  \int_a^b \Eu^{\iC N \Phi(s;\lambda)} \chi(s) \d s$.
  The integral $I(N;0)$ 
  reduces to the case treated in the first step
  with $k = \varphi''(s_*)/2$ and we have
  \begin{multline}
    \label{eq:36}
    I(N;0) = \Eu^{\iC N \varphi(s_*)} \chi(s_*) 
    \sqrt{\frac{2\pi}{N \varphi''(s_*)}} \Eu^{\iC \pi/4}
    + O\tonde{\frac{1+|b-a|}{N \delta}} = \\
    = \frac{J_* \Eu^{\iC N \varphi(s_*)}}{\sqrt{N}}
    + O\tonde{\1N}.
  \end{multline}
  In order to control $I(N;1)$,
  it will be enough to provide uniform bounds
  of order $O(1/N)$
  for the derivative $\de_\lambda I(N;\lambda)$.
  We derive~\eqref{eq:35} with respect to $\lambda$ 
  and then integrate by parts twice,
  \begin{multline}
    \label{eq:37}
    \de_\lambda I(N;\lambda) = \int
    \iC N \Phi_\lambda \Eu^{\iC N \Phi} \chi \d s =
    \int \iC N \Phi_s \Eu^{\iC N \Phi} \tonde{%
      \frac{\Phi_\lambda}{\Phi_s} \chi} \d s = \\
    = \int \tonde{\de_s \Eu^{\iC N \Phi}} \tonde{%
      \frac{\Phi_\lambda}{\Phi_s} \chi} \d s
    = -\int \Eu^{\iC N \Phi} \de_s \tonde{%
      \frac{\Phi_\lambda}{\Phi_s} \chi} \d s = \\
    = - \int \1{\iC N \Phi_s} \tonde{\de_s \Eu^{\iC N \Phi}}
    \de_s \tonde{\frac{\Phi_\lambda}{\Phi_s} \chi} \d s
    = \1{\iC N} \int \Eu^{\iC N \Phi} \de_s \tonde{%
      \1{\Phi_s} \de_s \tonde{%
        \frac{\Phi_\lambda}{\Phi_s} \chi}} \d s.
  \end{multline}
  There is no contribution coming from the boundary points,
  for we have chosen $\chi$ smooth and compactly supported
  in the interval $]a,b[$.
  Here we have used subscript notation
  for derivatives of $\Phi$,
  \begin{align*}
    \Phi_\lambda (s;\lambda) = \de_\lambda \Phi(s;\lambda) &=
    \varphi(s) - \varphi(s_*)
    - \frac{\varphi''(s_*)}{2} (s-s_*)^2
    = (s-s_*)^3 H(s), \\
    \Phi_s(s;\lambda) = \de_s \Phi(s;\lambda) &=
    \lambda \varphi'(s) + (1-\lambda) \varphi''(s_*)(s-s_*)
    =(s-s_*) K(s, \lambda),
  \end{align*}
  where
  \begin{align*}
    H(s) &= \int_0^1 \varphi'''\tonde{s_* + \theta (s-s_*)}
    \frac{(1-\theta)^2}{2} \d\theta, \\
    K(s;\lambda) &= (1-\lambda) \varphi''(s_*) + \lambda
    \int_0^1 \varphi''\tonde{s_* + \theta (s-s_*)} \d\theta.
  \end{align*}
  The condition $\varphi'' \ge 1$ implies $K \ge 1$.
  We have
  \begin{equation*}
    \frac{\Phi_\lambda}{\Phi_s} \chi = 
    (s-s_*)^2 \frac{H}{K} \chi =
    \tonde{\Phi_s}^2 \frac{H}{K^3} \chi = 
    \tonde{\Phi_s}^2 \widetilde{\chi},
  \end{equation*}
  The quantity
  \begin{equation*}
    \widetilde{\chi}(s;\lambda) =
    \frac{H(s)}{K^3(s;\lambda)} \chi(s)
  \end{equation*}
  is of class $C^2$ when $\varphi \in C^5$ and $\chi \in C^2$;
  moreover, we have uniform bounds
  for $\de_s^m \widetilde{\chi}$, $m=1,2$,
  in terms of uniform bounds of $\de_s^{k+3} \varphi$
  and $\de_s^k \chi$, with $k=0,1,2$.
  These computations show that the integrand 
  in the last integral of~\eqref{eq:37} is bounded,
  indeed we have
  \begin{equation*}
    \de_s \tonde{\1{\Phi_s} \de_s \tonde{%
        \frac{\Phi_\lambda}{\Phi_s} \chi}} = 
    \de_s \tonde{\1{\Phi_s} \de_s \tonde{%
        \tonde{\Phi_s}^2 \widetilde{\chi}}} = 
    \de_s \tonde{%
      2 \widetilde{\chi} + \Phi_s \widetilde{\chi}}.
  \end{equation*}
  It follows that $\de_\lambda I(N;\lambda) = O(1/N)$,
  uniformly with respect to $\lambda$.
  Hence, using~\eqref{eq:36} we obtain
  \begin{equation*}
    I(N) = I(N;1) = I(N,0)
    + \int_0^1 \de_\lambda I(N;\lambda) \d\lambda =
    \frac{J_* \Eu^{\iC N \varphi(s_*)}}{\sqrt{N}}
    + O\tonde{\1N}.
  \end{equation*}

  In the last step, we remove the assumption 
  on $\chi$ compactly supported in $]a,b[$
  by showing that the contribution from the boundary
  is also of order $O(1/N)$.
  Let $\beta \in C_0^\infty(]a,b[)$
  be such that $\beta(s)=1$ if $s \in [a+\delta, b-\delta]$
  and $\beta(s) = 0$ is $s \notin ]a+\delta/s, b-\delta/2[$.
  If we replace $\chi$ with $\beta\chi$ in the integral $I(N)$
  then we are in the situation of the second step.
  It remains to prove that the integral
  \begin{equation*}
    \widetilde{I}(N) =
    \int_a^b \Eu^{\iC N \varphi(s)} (1-\beta(s)) \chi(s) \d s
  \end{equation*}
  is of order $O(1/N)$.
  This follows immediately by an integration by parts,
  as in the proof of lemma~\ref{lemma:3},
  once we observe that when $s$ lies in support of $1-\beta$
  we have
  \begin{equation*}
    \abs{\varphi'(s)} =
    \abs{\int_{s_*}^s \varphi''(\sigma) \d\sigma} \ge
    \abs{s - s_*} \ge \frac{\delta}{2}.
  \end{equation*}
  Hence,
  \begin{equation*}
    \widetilde{I}(N) = \int_a^b \1{\iC N \varphi'}
    \tonde{\de_s \Eu^{\iC N \varphi}} (1-\beta) \chi \d s = 
    - \1{\iC N} \int_a^b \Eu^{\iC N \varphi}
    \de_s \tonde{\frac{(1-\beta)\chi}{\varphi'}} \d s,
  \end{equation*}
  which shows that
  \begin{equation*}
    \abs{\widetilde{I}(N)} \le \1N \cdot
    \int_a^b \quadre{%
      \frac{2 \de_s\tonde{(1-\beta)\chi}}{\delta}
      + \frac{4(1-\beta)\chi\varphi''}{\delta^2}} \d s.
  \end{equation*}
\end{proof}

\end{document}